\documentclass[12pt,a4paper,oneside,final,dvipdfm,dvips,pdfm]{amsart}

\usepackage{amssymb}
\usepackage{amsmath}
\usepackage{amscd}
\usepackage{graphicx}
\usepackage{graphics}
\usepackage{amsfonts}
\usepackage{pb-diagram}

\usepackage{amsthm,lscape}

\numberwithin{equation}{section}

\newtheorem{defn}{Definition}

\begin{document}

\title{Dynamical Systems and Topological Surgery}

\author{S. Antoniou}
\address{Department of Mathematics,
National Technical University of Athens,
Zografou campus, GR-157 80 Athens, Greece.}
\email{funkadelicgr@yahoo.gr}

\author{S. Lambropoulou}
\address{Department of Mathematics,
National Technical University of Athens,
Zografou campus, GR-157 80 Athens, Greece.}
\email{sofia@math.ntua.gr}
\urladdr{http://www.math.ntua.gr/$\tilde{~}$sofia}


\subjclass{57,37}

\date{}
\maketitle

\begin{abstract}
In this paper we try to establish a connection between a three-dimensional Lotka--Volterra dynamical system and two-dimensional topological surgery. There are many physical phenomena exhibiting  two-dimensional topological surgery through a `hole drilling' process. By our connection, such phenomena may be modelled mathematically by the above dynamical system.
\end{abstract}

\section{Introduction}

Topological surgery is a technique used for changing the homeomorphism type of a topological space. For example, all orientable surfaces may arise from the 2-dimensional sphere using surgery. Further, the phenomenon of surgery appears in many physical processes. In Section~2 we discuss in detail the mathematical definition of topological surgery in dimensions one, two and three and we give examples of physical processes exhibiting surgery in  dimensions one and two, such as: DNA recombination, reconnection of cosmic magnetic lines, cell mitosis, formation of tornadoes and whirls, falaco solitons, etc.

In Section~3 we present the three-dimensional Lotka--Volterra dynamical system, given and analyzed by N.Samardzija and D.Greller in \cite{SaGr1} (see system 3.1).  This system generalizes the classical two-dimensional Lotka--Volterra system \cite{Lo,Vo} and it depends on three parameters. As it turns out, in certain parametrical regions the system exhibits chaotic solutions shaping a fractal toroidal bead in 3-dimensional phase space. Moreover, with a continuous change of parameters, the chaotic attractor passes from spherical to toroidal shape via a `drilling' process along a slow manifold, thus performing 2-dimensional topological surgery. This connection was first observed in \cite{La} and consequently worked in \cite{An}. In this paper we establish this connection between the parametrical behavior of the solutions and the 2-dimensional topological surgery by solving the system numerically in each parametrical region and demonstrating this `hole drilling' process. This means that system 3.1 (more precisely, a one-parameter family) can model 2-dimensional surgery.

Further, there are physical phenomena resembling 2-dimensional surgery through a `drilling' process, such as formation of tornadoes and whirls, falaco solitons, etc. These applications are discussed in Section~4. By our connection, such phenomena can be modelled by system 3.1. Finally, in Section~5 we discuss the vectorial dynamo equation, the solutions of which we conjecture to model 3-dimensional surgery. This will be the subject of future work.

\section{The topological surgery}

\subsection{}   A {\it homeomorphism} between two $m$-dimensional topological spaces is a continuous bijection, such that the inverse is also a continuous map.  In Topology, two homeomorphic topological spaces are considered the same. For this reason they may be attached together and a homeomorphism between them can be used as `glue'.
\smallbreak

An  {\it $m$-manifold} $M$ is a metric space which may be covered by open sets, each of which is homeomorphic with ${\Bbb R}^m$ or the half--space ${\Bbb R}_+^m$ (in case $M$ has boundary). The aim of Topology is to classify all $m$-manifolds up to homeomorphism, that is, to characterize and list all non-homeomorphic $m$-manifolds, for each $m$.  Here we restrict our interest to orientable manifolds.

 The only connected 1-manifold is the circle $S^1$. The connected, orientable 2-manifolds without boundary are the 2-dimensional sphere $S^2$ (only the surface, without the inside), the torus (without the inside), the torus with two holes, with three holes, etc. In three dimensions the basic 3-manifold is the 3-dimensional sphere $S^3$, which may be viewed as the Euclidean space ${\Bbb R}^3$ with all points at infinity compressed into one single point. (Analogous views are true for spheres of all dimensions). Other examples of 3-manifolds are the lens spaces, homology 3-spheres, etc.
\smallbreak

In order to create new topological spaces out of known ones we use the following technique.

\begin{defn} \rm A {\it topological surgery} on an $m$-manifold $M$ is the topological procedure of creating a new $m$-manifold $M'$ by removing the product set $S^n\times D^{m-n}$ of the $n$-dimensional sphere with the $(m-n)$-dimensional ball, taking the closure of the remaining manifold with boundary and replacing $S^n\times D^{m-n}$ by $D^{n+1}\times S^{m-n-1}$, using a `gluing' homeomorphism along the common boundary $S^n\times S^{m-n-1}$, that is:
\[M' = (\overline{M\setminus S^n\times D^{m-n}}) \cup_{S^n\times S^{m-n-1}} D^{n+1}\times S^{m-n-1} \]
\end{defn}

An excellent reference on the subject is \cite{Ro}. In this paper we shall use the term `$m$-dimensional topological surgery' with no risk of confusion about the sub-dimension $n$. One-dimensional surgery, for example, means that we remove from one circle or from two circles two arcs $S^0\times D^1$ ($S^0$ is two points) and we reconnect the four boundary points $S^0\times S^0$ in a different way, obtaining in the end two or one circle depending on the type of reconnection. View Figure~1 (cf. p. 247 \cite{Ro}).

\bigbreak
\begin{figure}[h!]
\begin{center}
\includegraphics[width=6cm, bb=0 0 153 78]{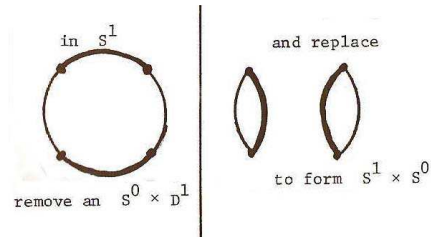}
\caption{Surgery in dimension 1}
{\footnotesize \underline{{\it Source}}: D. Rolfsen, {\em Knots and Links}, Publish or Perish Inc. (1976), p. 247.}
\end{center}
\end{figure}
\bigbreak

In this paper we are mainly interested in 2-dimensional surgery. From Definition~1 we have:
\[M' = (\overline{M\setminus S^0\times D^2}) \cup_{S^0\times S^1} D^1\times S^1 \]
That is, we remove two discs $S^0\times D^2$ from the surface $M$  and we replace them by a tube $D^1\times S^1$, which gets attached along the common boundary $S^0\times S^1$. This operation changes the homeomorphism type of the initial surface. In fact, according to the surface classification theorem, every connected, orientable 2-dimensional surface without boundary arises from the sphere $S^2$ by repeated surgeries. Note that, every time the above process is performed a new `handle' is added on the surface.

\smallbreak
For example, for $M=S^2$, the 2-dimensional sphere, the process is illustrated in Figure~2. Two points are specified first on the sphere, which are the centers of the discs to be removed. After the removal of the two discs $S^0\times D^2$ a tube $D^1\times S^1$ is attached along the common boundary $S^0\times S^1$, resulting in a torus. This is a different 2-dimensional surface, arising from the sphere via the surgery technique. Note that, if the tube were attached on the two-punctured sphere externally, the result would be a homeomorphic torus.  As another example, if the initial surface $M$ is two copies of the sphere $S^2$ and the two discs to be removed are one on each sphere, replacing the two discs by a tube yields one copy of the sphere $S^2$.

\begin{figure}[h!]
\begin{center}
\includegraphics[width=12cm, bb=0 0 170 42]{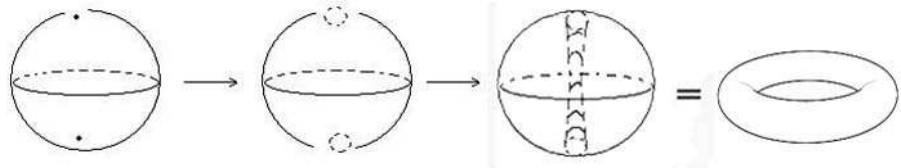}
\caption{Surgery on the sphere results in a torus}
\end{center}
\end{figure}

 Three-dimensional topological surgery is much harder to visualize than 2-dimensional surgery. Letting $m=3$ and $n=1$ in Definition~1, it means that we remove an embedded solid torus (a torus with the inside) $S^1\times D^2$  from a  3-manifold $M$ and we replace it by another solid torus $D^2\times S^1$ with the factors reversed. (To see that $S^1\times D^2$ is a solid torus consider a disc $D^2$ in every point of the circle $S^1$.) The replacing of the solid torus takes place by attaching along the common boundary $S^1\times S^1$ (which is a torus) via a gluing homeomorphism:
 \begin{equation} \label{3dsurgery}
M' = (\overline{M\setminus S^1\times D^2}) \cup_{S^1\times S^1} D^2\times S^1
\end{equation}
The result is a new 3-manifold $M'$.  Taking the 3-sphere $S^3 = {\Bbb R}^3 \cup \{\infty\}$ as our starting manifold, by known results in Topology, every 3-dimensional closed (compact, without boundary), connected, orientable manifold may arise as the result of surgery in  $S^3$ along a knot or link (these are embedded circles or copies of circles in three-space).

\subsection{}

Surgery is not just a technique used in Topology for changing the homeomorphism type of a topological space. It also happens in nature in diverse processes for ensuring new results. For example, 1-dimensional topological surgery happens in DNA recombination for changing the genetic sequence. View Figure~3 (cf. Wikipedia). It also happens in the magnetic reconnection, the phenomenon whereby cosmic magnetic field lines from different magnetic domains are spliced to one another, changing their patterns of connectivity with respect to the sources, see Figure~4 (cf. \cite{DaAn}).

\bigbreak
\begin{figure}[h!]
\begin{center}
\includegraphics[width=12cm, bb=0 0 269 29]{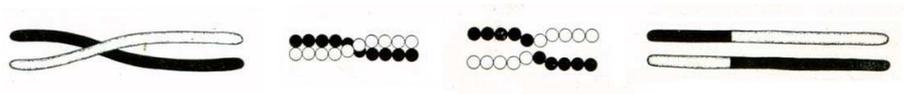}
\caption{Crossing over of chromosomes in DNA recombination}
{\footnotesize \underline{{\it Source}}: Wikipedia}
\end{center}
\end{figure}
\bigbreak

\bigbreak
\begin{figure}[h!]
\begin{center}
\includegraphics[width=12cm, bb=0 0 128 29]{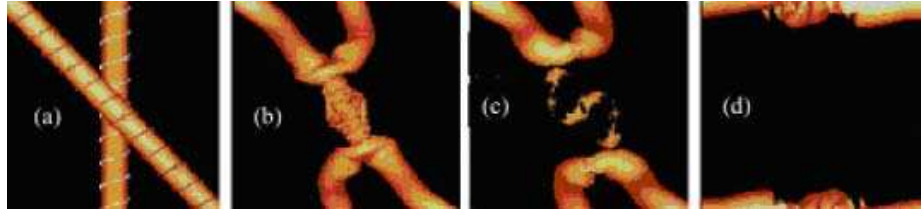}
\caption{The reconnection of cosmic magnetic lines}
{\footnotesize \underline{{\it Source}}: R.B. Dahlburg, S.K. Antiochos {\em Reconnection of Antiparallel Magnetic Flux Tubes}, J.  Geophysical Research {\bf 100}, No. A9 (1995) 16991--16998.}
\end{center}
\end{figure}
\bigbreak

Two-dimensional topological surgery happens, for example, in the biological process of mitosis, where a cell splits into two new cells. View Figure~5 (cf. p. 395 \cite{KeFa}). A similar picture can be found in the mechanism of gene transfer in bacteria (cf. p. 471 \cite{HHGRSV}). More examples are presented in Section~4.

\bigbreak
\begin{figure}[h!]
\begin{center}
\includegraphics[width=14cm, bb=0 0 612 192]{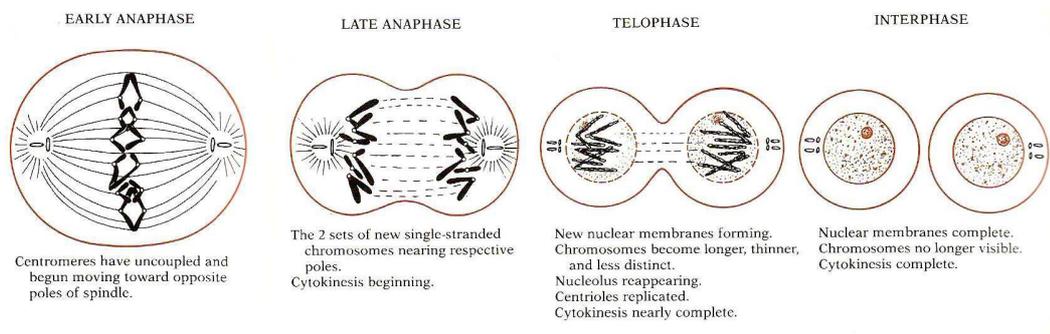}
\caption{The process of mitosis is an example of 2D surgery}
{\footnotesize \underline{{\it Source}}: W.T. Keeton, C.H. McFadden, {\em Elements of Biological Science}, W.W. Norton \& Company Inc., 3rd edition (1983), p. 395.}
\end{center}
\end{figure}
\bigbreak

Finally, a possible physical application of 3-dimensional topological surgery is discussed in Section~4.

\section{Connecting a 3--dimensional Lotka--Volterra model with 2-dimensional topological surgery}

\subsection{}
In \cite{SaGr1} Greller and Samardziya study the behavior of the following dynamical system that generalizes the classical Lotka--Volterra problem \cite{Lo,Vo} into three dimensions. We shall present the analysis of \cite{SaGr1} using the same notations and some of their figures. The system is a two-predator and one-prey model, where the predators $Y,Z$ do not interact directly with one another but compete for prey~$X$. Apart from a population model, the system may also serve as a biological model and a chemical model (cf. \cite{SaGr1}).

\begin{equation}
\left\{
\begin{array}{l}
\frac{dX}{dt}=X-XY+CX^2-AZX^2 \\
\\
\frac{dY}{dt}=-Y+XY \\
\\
\frac{dZ}{dt}=-BZ+AZX^2 \\
\end{array}
\right\} \ A, B, C \geq 0
\end{equation}

\vspace{.12in}

\noindent The parameters $A,B,C$ are analyzed in order to determine bifurcation properties of the system. In particular, the phenomenon of chaos evolving on a fractal torus is observed and analyzed by analytically derived two-dimensional stability diagrams. As parameters $A,B,C$ affect the dynamics of constituents  $X,Y,Z$, the authors were able to determine conditions for which the ecosystem of the three species results in steady, periodic or chaotic behavior. More precisely, the authors derive the following five steady state solutions for the system:

\begin{equation*}
S_{s1}=\left(
\begin{array}{c}
0 \\
0 \\
0
\end{array}
\right) ; \ \
S_{s2}=\left(
\begin{array}{c}
1 \\
1+C \\
0
\end{array}
\right) ; \ \
S_{s3}=\left(
\begin{array}{c}
\sqrt{B/A} \\
0 \\
\frac{1+C\sqrt{B/A}}{\sqrt{AB}}
\end{array}
\right) ;
\end{equation*}

\begin{equation*}
S_{s4}=\left(
\begin{array}{c}
-1/C \\
0 \\
0
\end{array}
\right) ; \ \
S_{s5}=\left(
\begin{array}{c}
-\sqrt{B/A} \\
0 \\
\frac{C\sqrt{B/A}-1}{\sqrt{AB}}
\end{array}
\right).
\end{equation*}

Since $X,Y,Z$ are populations, the negative steady states $S_{s4},S_{s5}$ are not taken into consideration. Figure~6     (cf. Fig.~1 in \cite{SaGr1}) illustrates the positions of steady states $S_{s1}, S_{s2}$ and $S_{s2}$ in $P^3$, where $P^3$ denotes the set of points in $\Bbb{R}^3$ with all coordinates $\geq0$.

\begin{figure}[h!]
\begin{center}
\includegraphics[width=10cm,bb=0 0 372 277]{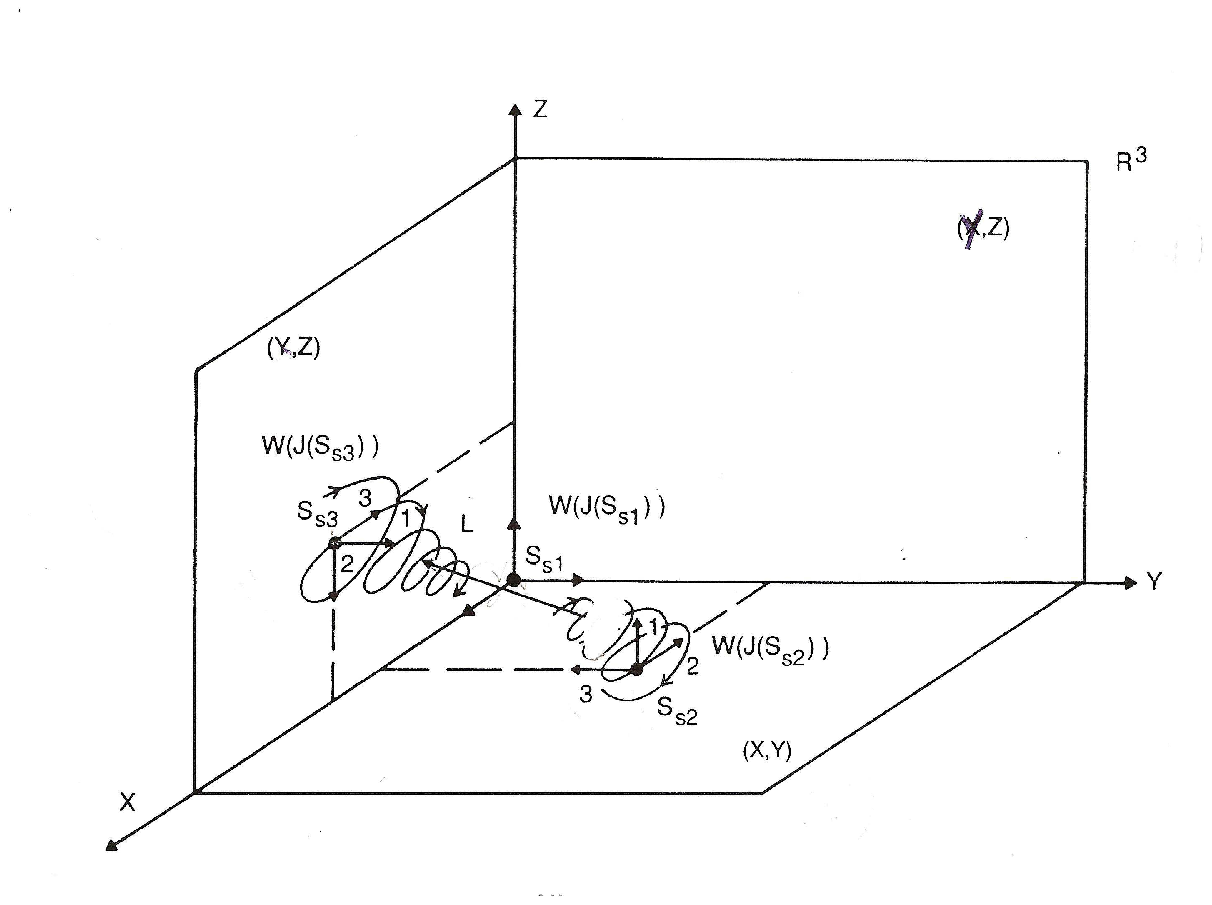}
\caption{Repelling steady state points $S_{s2}$ and $S_{s3}$.}
{\footnotesize \underline{{\it Source}}:  N. Samardzija, L. Greller, {\em Explosive route to chaos through a fractal torus in a generalized Lotka-Volterra Model}, Bull. Math. Biology {\bf 50}, No. 5 (1988), Fig.~1.}
\end{center}
\end{figure}

\smallbreak
Let, now, $J(S)$ be the Jacobian of system 3.1 evaluated at $S\in \Bbb{R}^3$ and let the sets $\Gamma\{J(S)\}$ and $W\{J(S)\}$ to be, respectively, the eigenvalues and corresponding associated eigenvectors of $J(S)$. Using the sets of eigenvalues and eigenvectors given below, the authors characterize the local behavior of the dynamical system around the points $S_{s1}$,$S_{s2}$, $S_{s3}$ using the Hartman-Grobman (or linearization) Theorem.

\begin{equation*}
\Gamma\{J(S_{s1})\}=\{1,-1,-B \}
 ; \ \
W\{J(S_{s1})\}=\left\{
\left[\begin{array}{c}
1 \\
0 \\
0
\end{array}\right],
\left[\begin{array}{c}
0 \\
1 \\
0
\end{array}\right],
\left[\begin{array}{c}
0 \\
0 \\
1
\end{array}\right]
\right\}
\end{equation*}

\begin{equation*}
\Gamma\{J(S_{s2})\}=\{A-B,(C+\sqrt{(C-2)^2-8})/2,(C-\sqrt{(C-2)^2-8})/2 \}
\end{equation*}

\begin{equation*}
W\{J(S_{s2})\}=\left\{
\left[\begin{array}{c}
1 \\
(C+1)/(A-B) \\
\frac{B+C-A+(C+1)/(B-A)}{A}
\end{array}\right],
\left[\begin{array}{c}
1 \\
\frac{C-\sqrt{(C-2)^2-8}}{2} \\
0
\end{array}\right],
\left[\begin{array}{c}
1 \\
\frac{C+\sqrt{(C-2)^2-8}}{2} \\
0
\end{array}\right]
\right\}
\end{equation*}

\vspace{.25in}
\noindent $\Gamma\{J(S_{s3})\}=$

\begin{equation*}
\ \ \ \ \ \ \ \ \ \ \ \ \left\{\sqrt{\frac{B}{A}}-1,\frac{-1+\sqrt{1-8B(1+C\sqrt{B/A})}}{2},\frac{-1-\sqrt{1-8B(1+C\sqrt{B/A})}}{2} \right\}
\end{equation*}

\vspace{.25in}
\noindent $W\{J(S_{s3})\}=$

\begin{equation*}
\left\{
\left[\begin{array}{c}
1 \\
 -1-\frac{2\sqrt{AB}(1+C\sqrt{B/A})}{\sqrt{B/A}-1}\\
\frac{2(1+C\sqrt{B/A})}{\sqrt{B/A}-1}
\end{array}\right],
\left[\begin{array}{c}
1 \\
0 \\
\frac{-1-\sqrt{1-8B(1+C\sqrt{B/A})}}{2B}
\end{array}\right],
\left[\begin{array}{c}
1 \\
0 \\
\frac{-1+\sqrt{1-8B(1+C\sqrt{B/A})}}{2B}
\end{array}\right]
\right\}
\end{equation*}

\vspace{.25in}
\bigbreak

Since $1>0$ and $-1,-B<0$, $S_{s1}$ is a saddle point for all $A,B,C$. By inspecting the variation of the eigenvalues for $S_{s2}$ and $S_{s3}$ in relation to $A,B,C$, the authors show that the behavior of $S_{s2}$ and $S_{s3}$ can be characterized by only two parameters, namely $C$ and $B/A$, as shown in the stability diagrams in Figure~7 (cf. Figs.~4, 5 and 6 in \cite{SaGr1}). Before examining the behavior of the system in each different parametrical region, the authors analyze it into three 2-dimensional subsystems, concluding that all solutions initiated in $P^3$ will remain in $P^3$ for all times. Furthermore, they show that system 3.1) is always bounded in $P^3$.


\begin{figure}[h]
\begin{center}
\includegraphics[width=8.5cm, bb=0 0 374 242]{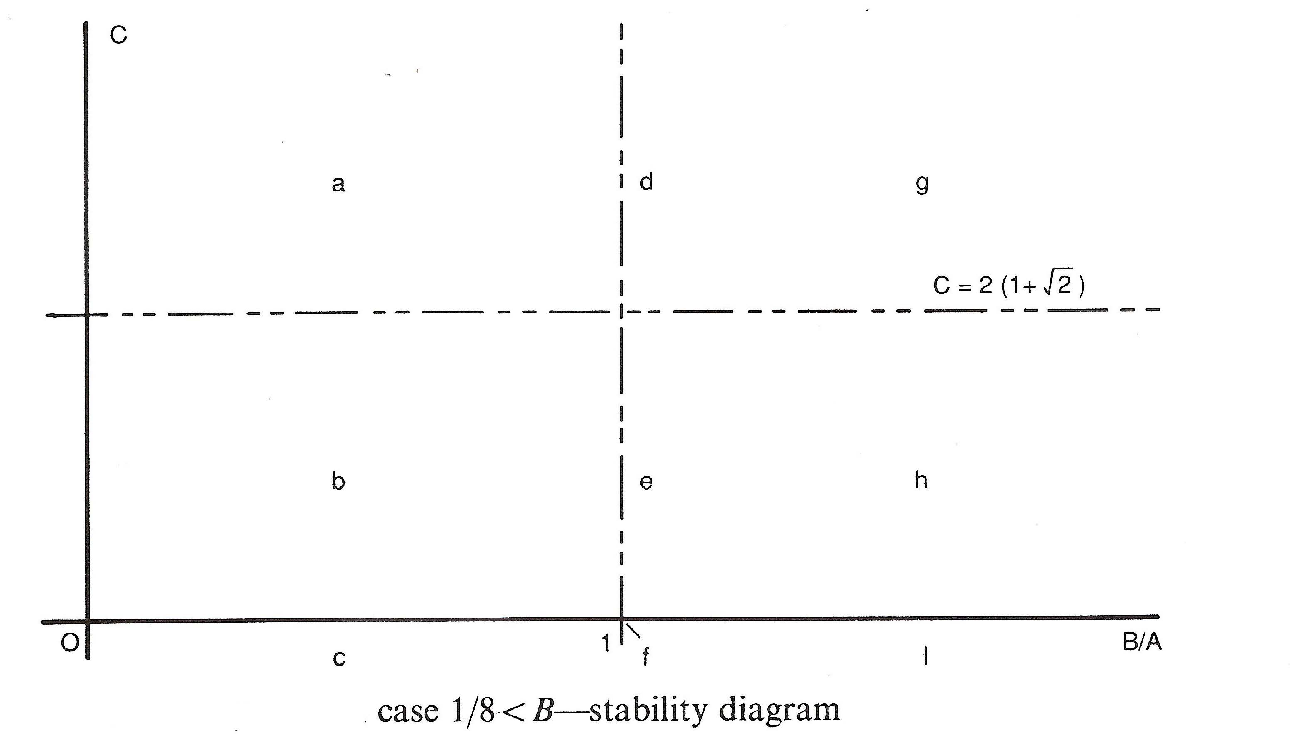}
\includegraphics[width=8.5cm, bb=0 0 354 255]{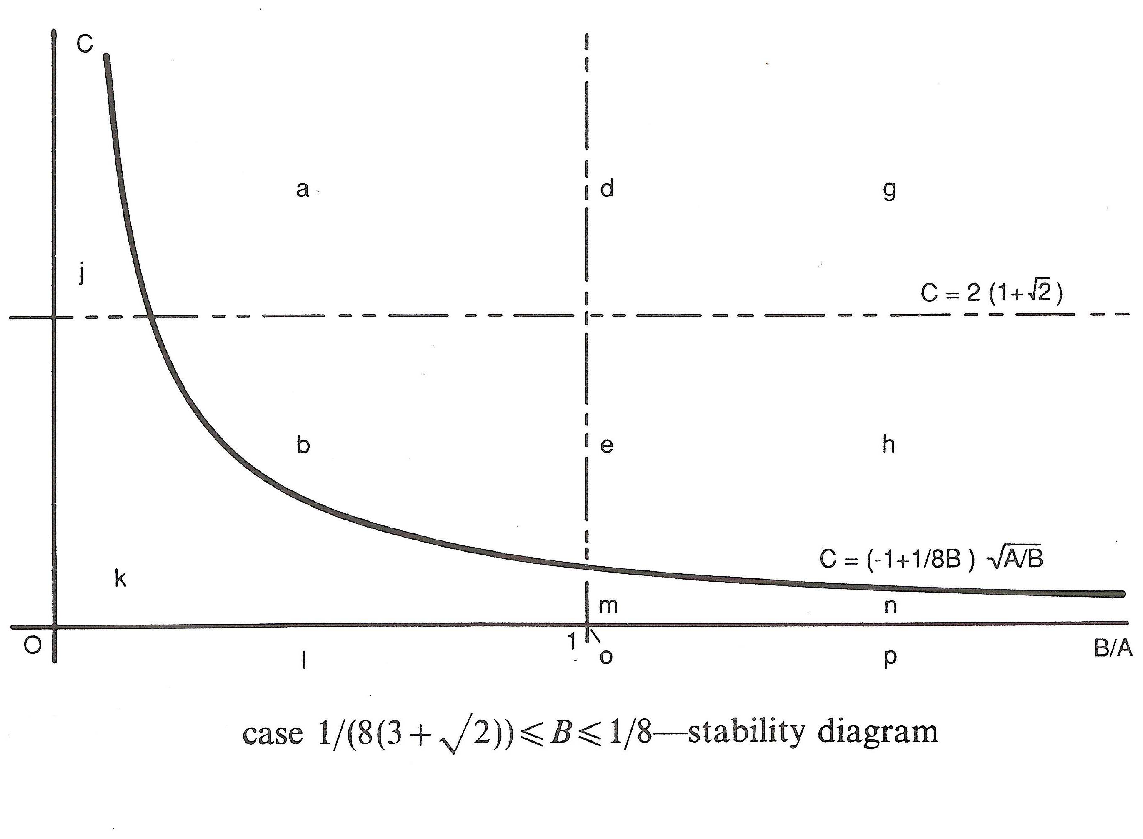}
\includegraphics[width=9cm, bb=0 0 373 237]{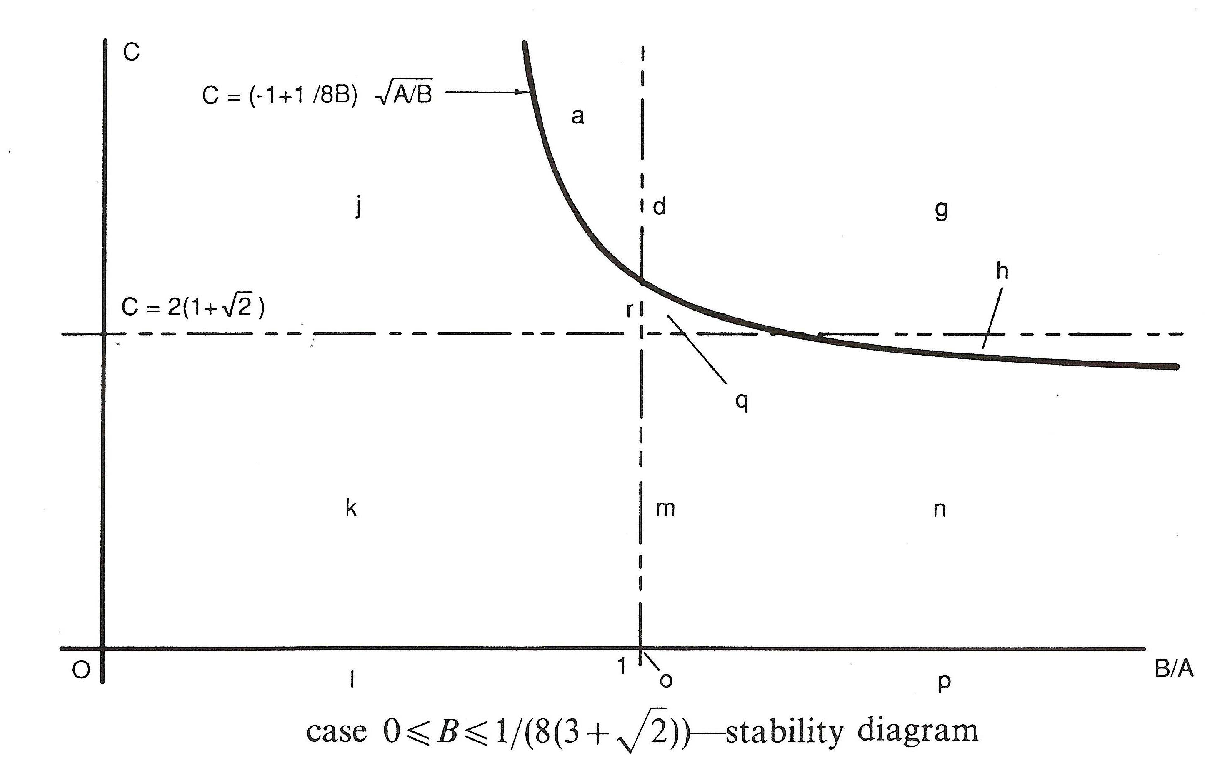}
\caption{Stability diagrams}
{\footnotesize \underline{{\it Source}}:  N. Samardzija, L. Greller, {\em Explosive route to chaos through a fractal torus in a generalized Lotka-Volterra Model}, Bull. Math. Biology {\bf 50}, No. 5 (1988), Figs.~4, 5 and 6.}
\end{center}
\end{figure}

Parametrical regions $a, b, c, d, e, f, j, k, l, m, o$ and $r$ are stable steady state regions, while $i$ and $p$ are periodic/stable steady state regions. Interesting behavior arises in the parametrical regions $g, h, n$ and $q$ where $S_{s2}$  and $S_{s3}$ are an inward instable vortex  and an outward stable vortex, respectively (view Figure~6). Let $\lambda_1, \lambda_2, \lambda_3$ denote the eigenvalues of either $S_{s2}$ or $S_{s3}$. In both cases they must satisfy the conditions $\lambda_1\in \Bbb{R}$ and $\lambda_2, \lambda_3\in \Bbb{C}$ with $\lambda_3= \lambda_2^*$, the conjugate of $\lambda_2$. The eigenvalues of $S_{s2}$ must further satisfy $\lambda_1<0$ and  $Re(\lambda_2)=Re(\lambda_3)>0$, while the eigenvalues of $S_{s3}$ must further satisfy $\lambda_1>0$ and  $Re(\lambda_2)=Re(\lambda_3)<0$.

Note now that the point $S_{s2}$ as well as the  eigenvectors corresponding to its two complex eigenvalues, all lie in the $xy$--plane. On the other hand, the point $S_{s3}$ and also the  eigenvectors corresponding to its two complex eigenvalues all lie in the $xz$--plane. Finally, the eigenvector corresponding to the real eigenvalue of $S_{s3}$ breaks out of the $xz$--plane  and redirects the flow towards $S_{s2}$. Thus, since we have two repelling steady points and since there are no sinks in $P^3$ and the system is bounded, the formation of chaotic/periodic type trajectories is expected in regions $g, h, n$ and $q$.

\bigbreak

Setting $B/A=1$ and equating the right side of system 3.1 to zero, one finds as solution the one-dimensional manifold
\[L=\{(X,Y,Z); \ X=1, Y+AZ=1+C \}\]
that passes through  $S_{s2}$  and $S_{s3}$. Since all points on $L$ are singular, there is no motion along it. If one moves in parameter space to the left of $B/A=1$, the manifold $L$ exerts flow toward the inward stable vortex $S_{s3}$ and if one moves to the right of  $B/A=1$, the manifold $L$ exerts motion toward the inward unstable vortex $S_{s2}$. Consequently, stable solutions are generated left of and including the line $B/A=1$, while the chaotic/periodic regions $g, h, n$ and $q$ that we are interested in appear on the right of the line $B/A=1$.

\bigbreak

Using CRAY FORTRAN, the authors solved system 3.1 numerically in chaotic/periodic regions $g, h, n$ and $q$ (where $B/A>1$) and showed that the trajectories are quickly trapped into an attractor between $S_{s2}$ and  $S_{s3}$. Since the inter-arrival times exhibit stochastic behavior, the authors state that the attractor is non-periodic, thus being chaotic. The line $B/A=1$ exhibits the slowest dynamic on the attractor, so $L$ is  referred to as the `slow manifold', while the `fast manifold' is the outer shell of the attractor. As parameter $A$ changes, the slow manifold $L$ opens into a hole and the chaotic attractor evolves into a conventional toroidal shape.

\subsection{}
We are now ready to make the connection of the two areas.
Using MATLAB ODE Solver we solved system 3.1 numerically and we represented the attractor in the chaotic/periodic regions $g, h, n$ and $q$. We too observed that, as parameter $A$ is increased, the slow manifold $L$ opens into a hole. Figure~8 illustrates instances of this process in chaotic region $q$ for $B=0.0145, C=5.5$ and for $A=0.01305, 0.01335, 0.01365, 0.01395$ and $0.01425$. Once the trajectory is trapped into the attractor, the outer shell and the line $L$ (slow manifold) are drawn.  In the first instance the trajectory is creating the 2-sphere and the slow manifold $L$, which links the two steady state  points $S_{s2}$  and $S_{s3}$. The line $L$ `selects' the centers of the two discs  $S^0\times D^2$ that will be removed through the surgery process. By increasing parameter $A$, the hole is drilled and our attractor passes from a spherical to a toroidal shape, so 2-dimensional surgery is performed. The same shape transition also happens in regions $g, h$ and $n$ for appropriate choices of parameters. We chose chaotic region $q$ because it demonstrates our point more clearly. Therefore, system 3.1 may perform 2-dimensional surgery in the above parametrical regions.

\begin{figure}[h!]
\begin{center}
\includegraphics[width=14cm,bb=0 0 363 475]{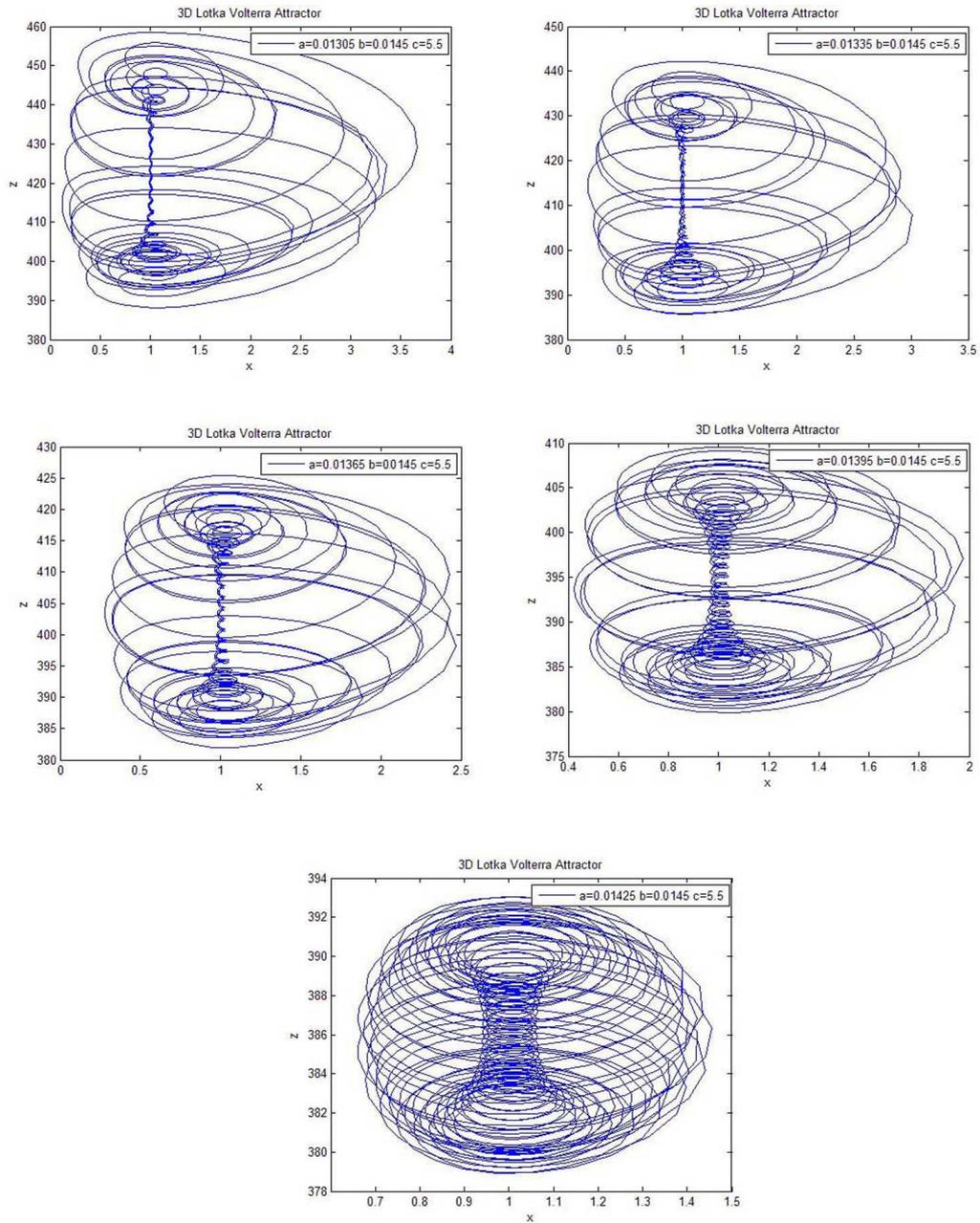}
\caption{Opening of the slow manifold $L$ into a hole as parameter $A$ increases.}
\end{center}
\end{figure}

\smallbreak
The use of different numerical methods may affect the shape of the attractor but the `hole drilling' process is a always a common feature.

\smallbreak
\noindent {\it Conclusion:} The (one-parameter family of the) three dimensional Lotka--Volterra system 3.1 models the process of 2-dimensional topological surgery. On the other hand, instead of viewing surgery as an abstract topological process, we may now view it as a property of an one-parameter dynamical system.

\section{Applications}

Many physical phenomena start by selecting two points (poles) and by joining them through a helicoid line. This process resembles surgery, regardless of the topological surface the two points belong to. Moreover, `drilling' with coiling is naturally chosen in various physical processes, since it is, from the energy point of view, the most economical way to open a hole. Thus, as a consequence of our connection, every physical phenomenon exhibiting the above process can be modelled by the three-dimensional Lotka--Volterra system 3.1.

\subsection{} A phenomenon exhibiting the above process is the creation of Falaco solitons. These are singular surfaces connected by means of a stabilizing invisible thread or string. For a full exposition see \cite{Ki}. In Figure~9 (cf. pp. 147, 150 \cite{Ki}) one can observe three pairs of these topological structures in a swimming pool.

\begin{figure}[h!]
\begin{center}
\includegraphics[width=10cm,bb=0 0 142 59]{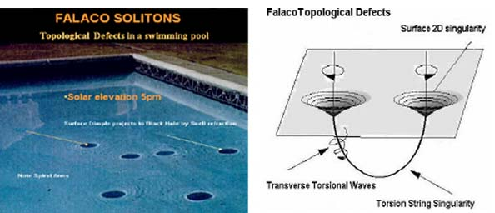}
\caption{Falaco Solitons in a swimming pool}
{\footnotesize \underline{{\it Source}}: R.M. Kiehn, {\em Non-Equilibrium Systems and Irreversible Processes}, Adventures in Applied  Topology {\bf 1}, Non Equilibrium Thermodynamics, University of Houston Copyright CSDC. INC, (2004), pp. 147, 150.}
\end{center}
\end{figure}

The invisible thread connecting the two singular surfaces (poles) of Falaco solitons corresponds to our slow manifold $L$. The hole drilling process seen in Figure~9 is exactly analogous to the one observed in Section~3 (recall Figure~8). Starting by the two poles, and by coiling the joining line, the hole becomes bigger and bigger until surgery is completed. In Figure~9 one can see that each soliton pair contains one inward instable vortex and one outward stable vortex, just like our two steady points $S_{s2}$  and $S_{s3}$ of system 3.1 (recall Figure~6).

\smallbreak
As mentioned in \cite{Ki}, p. 152 ``a Falaco soliton is a topological object that appears at all scales, so the concept of Falaco solitons is a universal phenomenon. At the microscopic level, the method offers a view of forming spin pairs that is different from Cooper pairs and could offer insight to superconductivity. At the level of cosmology, the concept of Falaco solitons can be found in solutions to the Navier--Stokes equations in a rotating frame of reference.  Based on the experimental creation of Falaco solitons in a swimming pool, it has been conjectured that M31 and the Milky Way galaxies could be connected by a topological thread". View Figure~10.

\begin{figure}[h!]
\begin{center}
\includegraphics[width=5cm,bb=0 0 85 72]{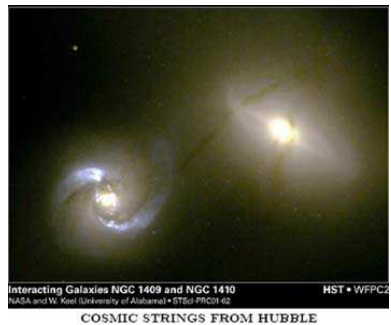}
\caption{The galaxies M31 and Milky Way}
\end{center}
\end{figure}

\subsection{} Another physical phenomenon which can be modelled by system 3.1 is the formation of tornadoes. As Figure~11  demonstrates, if certain meteorological conditions are met, funnel-shaped clouds start descending toward the ground. Once they reach it, they become tornados. This procedure can be viewed as a two-dimensional surgery. First the poles are chosen, one on the tip of the cloud and the other on the ground. Then, starting from the first point, the wind revolves in a  helicoidal motion toward the second point, until the hole is drilled and, thus, 2-dimensional surgery is performed. Note that, if one deforms the water surface of falaco solitons to a sphere, the visual correspondence to tornado structures between the earth surface and the tropopause is remarkable.

\begin{figure}[h!]
\begin{center}
\includegraphics[width=8cm,bb=0 0 1278 445]{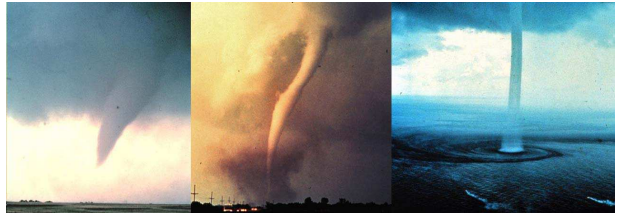}
\caption{Funnel clouds drilling and tornado formation}
\end{center}
\end{figure}

\subsection{} Another physical phenomenon completely analogous to the formation of tornadoes is the formation of whirls in the sea.

\subsection{} Finally, a phenomena which is very similar to the previous ones is a wormhole.
In physics, a wormhole is a hypothetical topological feature of space-time that is basically a `shortcut' through space and time. It has at least two mouths which are connected to a single throat or tube. View Figure~12.
 This structure was presented by Wheeler in 1955 but was considered to be unattainable in a practical sense.  While there is no observational evidence for wormholes, space-times containing wormholes are known to be valid solutions in general relativity.

\begin{figure}[h!]
\begin{center}
\includegraphics[width=6cm,bb=0 0 450 294]{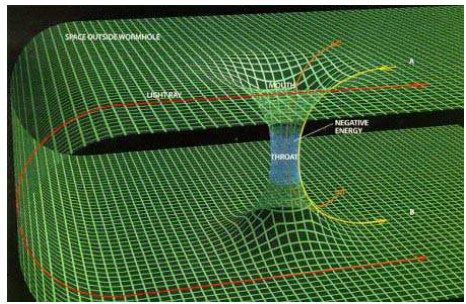}
\caption{A representation of a wormhole}
\end{center}
\end{figure}

\section{3-dimensional topological surgery and dynamo equations}

Recall now (\ref{3dsurgery}) that describes 3-dimensional topological surgery. Letting $M=S^3$, in order to cut out a solid torus from $S^3$ we need to know the following property: the complement of a solid torus $V_1$ in $S^3$ is another solid torus $V_2$, that contains the point at infinity, and one can write $S^3 = V_1\cup_\vartheta V_2$. The Hopf fibration of $S^3$, illustrated in Figure~13, gives an idea of this splitting of $S^3$. The surgery process starts and ends with two 3-dimensional manifolds than can be represented by two tori. The parallels of the torus of  $S^3$  become the new meridians of the torus  of our new 3-manifold.

\begin{figure}[h!]
\begin{center}
\includegraphics[width=8cm,bb=0 0 402 258]{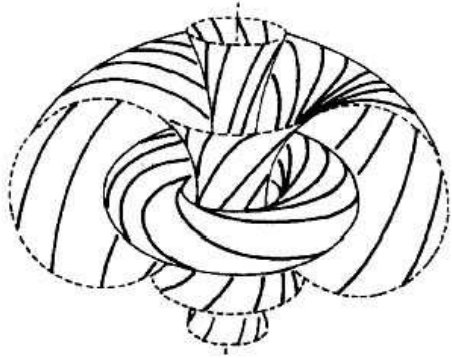}
\caption{Representation of $S^3$ as a union of two solid tori.}
\end{center}
\end{figure}

\smallbreak
On the other hand, similar toroidal structure can be found in dynamo theory, see for example \cite{Fi}. More precisely, when conducting fluid flows across an existing magnetic field, electric currents are induced which, in turn, create another magnetic field. When this new magnetic field reinforces the original magnetic field a self-sustained dynamo is created. This is called the ``Dynamo Theory" and it explains how the magnetic field of dense objects, such as the white dwarves, the sun or the earth, is sustained.


In order to simulate the Earth's magnetic field, Glatzmaier and Roberts solved the dynamo equation numerically \cite{GlRo1}. Their result is precisely a three dimensional manifold represented by two tori, view Figure~14.

\begin{figure}[h!]
\begin{center}
\includegraphics[width=5cm,bb=0 0 171 188]{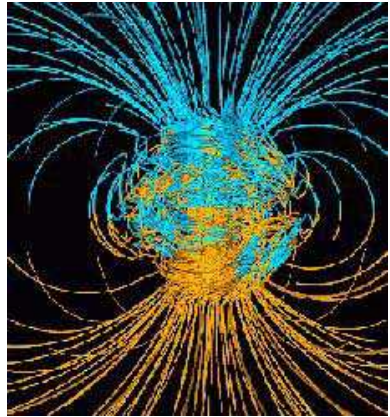}
\caption{Numerical solution of the dynamo equation}
{\footnotesize \underline{{\it Source}}:  G.A. Glatzmaier, P.H. Roberts, {\em Three-dimensional self-consistent computer simulation of the geomagnetic field reversal}, Nature {\bf 377} (1995).}
\end{center}
\end{figure}

\smallbreak
We believe that dynamo theory is the key to our next step: modelling 3-dimensional surgery. This will be the subject of future work.


\end{document}